\documentclass[12pt,a4paper]{article}

\usepackage[latin1]{inputenc}
\usepackage{latexsym,amssymb,stmaryrd}
\usepackage{amsmath,amssymb,amsfonts,amscd,theorem}
\usepackage{enumerate,hyperref}
\usepackage[francais,english]{babel}

\newcommand{\Ker}{\operatorname{Ker}}

\renewcommand{\hom}{\operatorname{Hom}}

\newcommand{\CC}{\mathbb{C}}

\newcommand{\image}{\operatorname{Im}}
\newcommand{\rk}{\operatorname{rk}}
\newcommand{\tr}{\operatorname{tr}}

\newcommand{\QQ}{\mathbb{Q}}

\newcommand{\ZZ}{\mathbb{Z}}
\newcommand{\Ad}{\operatorname{Ad}}
\newcommand{\Aut}{\operatorname{Aut}}

\def\co{\colon\thinspace}
\newenvironment{proof}[1][Proof]%
{
\begin{trivlist} \item[]  {\em #1.} }%
{\hspace*{\fill} $\Box$
\end{trivlist}}

{\bigl(
\begin{smallmatrix} }%
{
\end{smallmatrix} \bigr)}

\newtheorem{thm}{Theorem}
\newtheorem{prop}[thm]{Proposition}
\newtheorem{lemma}[thm]{Lemma}

\newtheorem{cor}[thm]{Corollary}
\newtheorem{rem}[thm]{Remark}

\theorembodyfont{\rmfamily}
       \newtheorem{definition}[thm]{Definition}
       \newtheorem{remark}[thm]{Remark}

\title{Deformations of metabelian representations of  knot groups into $SL(3,\CC)$}

\author{Leila Ben Abdelghani, Michael Heusener and Hajer Jebali}
\date{}

\begin{document}
\maketitle

\begin{abstract}
Let $K$ be a knot in $S^3$ and $X$ its complement. We study deformations of reducible metabelian representations of the knot group $\pi_1(X)$ into $SL(3,\CC)$ which are associated to a double root of the Alexander polynomial. We prove that these reducible metabelian representations are smooth points of the representation variety and that they have irreducible non metabelian deformations.

\end{abstract}
\section*{Introduction}
\label{sec:intro}
Let $K$ be a knot in $S^3$. We let  
$X=\overline{ S^3\smallsetminus V(K)}$ denote the knot complement where $V(K)$ is a tubular neighbourhood of $K$.
Moreover we let $\pi=\pi_1(X)$ denote the fundamental group of $X$. Let $\mu$ be a meridian of $K$ and let $\Delta_K(t)\in\ZZ[t,t^{-1}]$ denote the Alexander polynomial of $K$. We associate to $\alpha\in\CC^*$ a homomorphism \begin{align*}
\eta_\alpha\co &\pi\to\CC^*\\
&\gamma\mapsto\alpha^{|\gamma|}
\end{align*}
with $|\gamma|=p(\gamma)$, where $p\co\pi\to\pi/\pi'\simeq\ZZ$ denotes the canonical projection. Note that $\eta_\alpha(\mu)=\alpha$. We define $\CC_\alpha$ to be the $\pi$-module $\CC$ with the action induced by $\alpha$, i.e.\ $\gamma\cdot x=\alpha^{|\gamma|}x$, for all $x\in\CC$ and for all $\gamma\in\pi$. The trivial $\pi$-module $\CC_1$ is simply denoted $\CC$.

Burde and de Rham proved, independently, that when $\alpha$ is a root of the Alexander polynomial there exists a reducible metabelian, non abelian, representation $\phi\co\pi\to GL(2,\CC)$ such that \[\phi(\gamma)=\begin{pmatrix}
\alpha^{|\gamma|}&z(\gamma)\\
0&1
\end{pmatrix}\,.\]
Here $z$ is a $1$-cocycle in $Z^1(\pi,\CC_\alpha)$ which is not a coboundary (see \cite{Bur} and \cite{dR}). The homomorphism $\phi$ induces a representation into $SL(2,\CC)$ (Lemma~\ref{lem:homo}) given by
\[\tilde{\phi}(\gamma)=\alpha^{-1/2}(\gamma)\phi(\gamma),\quad\forall\gamma\in\pi\,,\]
where $\alpha^{-1/2}\co\pi\to\CC^*$ is a homomorphism such that $(\alpha^{-1/2}(\gamma))^2=\alpha^{-|\gamma|}$, for all $\gamma\in\pi$.

This constitutes the starting point to the study of the problem of deformations of metabelian and abelian representations in $SL(2,\CC)$ or $SU(2)$ that correspond to a simple zero of the Alexander polynomial (see \cite{F-K} and \cite{H-P-S}). The result of \cite{F-K} is generalized in \cite{H} and \cite{H-K} by replacing the condition of the simple zero by a condition on the signature operator. Similar results are established in \cite{S} in the case of cyclic torsion, but unfortunately none of these results has been published. In \cite{BA-thesis}, \cite{BA} and \cite{BA-L}, the authors considered the case of any complex connected reductive or real connected compact Lie group. They supposed that the $(t-\alpha)$-torsion of the Alexander module is semisimple. A $PSL(2,\CC)$ version was given recently in \cite{H-P}.

Throughout this paper, we suppose that $\alpha\in\CC^*$ is a multiple root of the Alexander polynomial $\Delta_K(t)$ and that 
$\dim H^1(\pi,\CC_\alpha)=1$. This means that the $(t-\alpha)$-torsion of the Alexander module is cyclic of the form \[\tau_\alpha=\CC[t,t^{-1}]\big/(t-\alpha)^r \text{ where } r\geq 2 \,.\]
As particular generalization of Burde and de Rham's result it is established in \cite{Hajer} that in this case
there exists a reducible metabelian, non abelian, representation $\rho_0\co\pi\to GL(3,\CC)$ defined by
\[\rho_0(\gamma)=\begin{pmatrix}
\alpha^{|\gamma|}&z(\gamma)&g(\gamma)\\
0&1&h(\gamma)\\
0&0&1
\end{pmatrix}\,.\]
Here $h\co\pi\to(\CC,+)$ is a non trivial homomorphism and 
$g\co\pi\to\CC_\alpha$ is a $1$-cochain verifying $\delta g+z\cup h=0$. We normalize $\rho_0$ by considering
\begin{align*}
\tilde{\rho}\co\pi & \to SL(3,\CC)\\
\gamma &\mapsto\alpha^{-1/3}(\gamma)\rho_0(\gamma)
\end{align*}
where $\alpha^{-1/3}\co\pi\to\CC^*$ is a homomorphism (Lemma~\ref{lem:homo}) such that 
\[\left(\alpha^{-1/3}(\gamma)\right)^3=\alpha^{-|\gamma|},\quad\text{for all}\ \gamma\ \in\pi\,.\]

The Lie algebra $sl(3,\CC)$ turns into a $\pi$-module via the adjoint action of the representation $\tilde{\rho}$, $\Ad\circ \tilde\rho\co \pi\to\Aut(sl(3,\CC))$. The aim of this paper is to answer the following question: when can $\tilde{\rho}$ be deformed into irreducible non metabelian representations?

We use the technical approach of \cite{H-P} to prove the following result:
\begin{thm}\label{mainresult}
We suppose that the $(t-\alpha)$-torsion of the Alexander module is cyclic of the form $\CC[t,t^{-1}]\big/(t-\alpha)^r$.
If $\alpha$ is a double root of the Alexander polynomial i.e.\ $r=2$, then there exist irreducible non metabelian representations from $\pi$ into $SL(3,\CC)$ which deform $\tilde{\rho}$. Moreover, the representation $\tilde{\rho}$ is a smooth point of the representation variety $R(\pi,SL(3,\CC))$; it is contained in an unique $10$-dimensional component $R_{\tilde\rho}$ of $R(\pi,SL(3,\CC))$.
\end{thm}

The first example of classical knots whose Alexander polynomial has a double root $\alpha$ such that the $(t-\alpha)$-torsion of the Alexander module is cyclic is $8_{20}$.

This paper is organized as follows: In Section~\ref{notations} the basic notation and facts are presented. The Section~\ref{deform} includes the proof of Theorem~\ref{mainresult}. The cohomology computations are done in Section~\ref{cohometab}. The aim of Section~\ref{nature} is to study the nature of the deformations of $\tilde\rho$. 

\paragraph{Acknowledgements.} The authors are pleased to acknowledge the support by the French-Tunisian CMCU project \no 06S/1502.  
Moreover, the authors like to thank the anonymous referee for her/his valuable remarks. Her/his suggestions helped to improve in a substantial way the exposition and the content of the paper.

\section{Notation and facts}
\label{notations}

\begin{lemma}\label{lem:homo}
Let $n\geq 2$ and let $\eta\co\pi\to\CC^*$ be  a homomorphism. Then there exists $\tilde{\eta}\co\pi\to\CC^*$ a homomorphism such that 
$\left(\tilde{\eta}(\gamma)\right)^n=\eta(\gamma),\ \text{for all}\ \gamma\in\pi$.
\end{lemma}

\begin{proof}
Let $\lambda\co\pi\to\CC^*$ be a map satisfying $\left(\lambda(\gamma)\right)^n=\eta(\gamma)$, for all $\gamma\in\pi$. Then there exists a map $\omega\co\pi\times\pi\to U_n$ such that
\[\lambda(\gamma_1\gamma_2)=\lambda(\gamma_1)\lambda(\gamma_2)\omega(\gamma_1,\gamma_2)\quad,\ \forall\gamma_1,\gamma_2\in\pi\,,\]
where $U_n=\{\xi\in\CC^*\ |\ \xi^n=1\}\simeq\ZZ/n\ZZ$. Hence $\lambda$ is unique up to multiplication by a $n$-th root of unity.
It is not hard to check that $\omega$ is a $2$-cocycle in $Z^2(\pi,\ZZ/n\ZZ)$. Since $H^2(\pi,\ZZ/n\ZZ)=0$, there exists a $1$-cochain $d$ such that $\omega=\delta d$ and we can easily verify that $\lambda d$ is a homomorphism satisfying
\[(\lambda(\gamma) d(\gamma))^n=\eta(\gamma),\quad\forall\gamma\in\pi\,.\]
\end{proof}

\subsection{Group cohomology}
\label{prod}

The general reference for this section is Brown's book \cite{Br}.
Let $A$ be a $\pi$-module. We denote by $C^*(\pi,A)$ the cochain complex. An element of $C^n(\pi,A)$ can be viewed as a function $f\co\pi^n\to A$, i.e. as a function of $n$ variables from $\pi$ to $A$. The coboundary operator $\delta\co C^n(\pi,A)\to C^{n+1}(\pi,A)$ is given by:
\begin{align*}\delta f(\gamma_1,\ldots,\gamma_{n+1})&=\gamma_1\cdot f(\gamma_2,\ldots,\gamma_{n+1})+\displaystyle\sum_{i=1}^{n}(-1)^i f(\gamma_1,\ldots,\gamma_{i-1},\gamma_i\gamma_{i+1},\ldots,\gamma_{n+1})\\
&+(-1)^{n+1} f(\gamma_1,\ldots,\gamma_n)\,.\end{align*}
Note that $C^0(\pi,A)\simeq A$ and that, for $a\in C^0(\pi,A)$, we have:
\[\delta a(\gamma)=(\gamma-1)\cdot a\,,\quad\forall\ \gamma\in\pi\,.\]
 The coboundaries (respectively cocycles, cohomology) of $\pi$ with c\oe fficients in $A$ are denoted by $B^*(\pi,A)$ (respectively $Z^*(\pi,A),\ H^*(\pi,A)$). For $z$ a cocycle in $Z^i(\pi,A)$, $i\geq 1$, the cohomology class in $H^i(\pi,A)$ is denoted $\{z\}$.

Let $A_1,\ A_2$ and $A_3$ be $\pi$-modules. The cup product of two cochains $u\in C^p(\pi,A_1)$ and $v\in C^q(\pi,A_2)$ is the cochain $u\cup v\in C^{p+q}(\pi,\ A_1\otimes A_2)$ defined by
\[u\cup v(\gamma_1,\ldots,\gamma_{p+q}):=u(\gamma_1,\ldots,\gamma_p)\otimes\gamma_1\ldots\gamma_p\circ v(\gamma_{p+1},\ldots,\gamma_{p+q})\,.\]

Here $A_1\otimes A_2$ is a $\pi$-module via the diagonal action. It is possible to combine the cup product with any bilinear map $A_1\otimes A_2\to A_3$. We are only interested by the product map $\CC\otimes\CC_{\alpha^{\pm}}\to\CC_{\alpha^\pm}$ and $\CC_\alpha\otimes\CC_{\alpha^{-1}}\to\CC$.

A short exact sequence
\[ 0\to A_1 \stackrel{i}{\longrightarrow} A_2 \stackrel{p}{\longrightarrow} A_3 \to 0\]
of $\pi$-modules gives rise to a short exact sequence of cochain complexes:
\[ 0 \to C^*(\pi,A_1)\stackrel{i_*}{\longrightarrow} C^*(\pi,A_2) \stackrel{p_*}{\longrightarrow} C^*(\pi,A_3) \to 0\,.\]
In the sequel we will make use of the corresponding long exact cohomology sequence:
\[ 0\to H^0(\pi,A_1)\longrightarrow H^0(\pi,A_2) \longrightarrow H^0(\pi,A_3)\stackrel{\delta^1}{\longrightarrow}
H^1(\pi,A_1)\longrightarrow \cdots\,.\]
In order to define the connecting homomorphism 
$\delta^{n+1}\co H^n(\pi,A_3)\longrightarrow
H^{n+1}(\pi,A_1)$ we let $\delta_2$ denote the coboundary operator of $C^*(\pi,A_2)$. If $z\in Z^n(\pi,A_3)$ is a cocycle then $\delta^{n+1}(\{z\}) = \{ i_*^{-1} ( \delta_2 (\tilde z))\}$ where
the cochain $\tilde z \in p_*^{-1}(z)\subset C^n(\pi,A_2)$ is any lift of $z$.

\subsection{Group cohomology and representation varieties}
\label{cohorep}
The set $R_n(\pi):=R(\pi,SL(n,\CC))$ of homomorphisms of $\pi$ in $SL(n,\CC)$ is called the representation variety of $\pi$ in $SL(n,\CC)$ and is a (not necessarily irreducible) algebraic variety. 
\begin{definition}
A representation $\rho\co\pi\to SL(n,\CC)$ of the knot group $\pi$ is called abelian (resp. metabelian) if the restriction of $\rho$ to the first (resp. second) commutator subgroup of $\pi$, denoted $\pi'$ (resp. $\pi''$) is trivial.
\end{definition}

In this section, we present some results of \cite{H-P} that we will use in the sequel.
Let $\rho\co\pi\to SL(n,\CC)$ be a representation. 
The Lie algebra $sl(n,\CC)$ turns into a $\pi$-module via 
$Ad\circ\rho$. This module will be simply denoted by 
$sl(n,\CC)_\rho$. A cocycle $d\in Z^1(\pi,sl(n,\CC)_\rho)$ is a map 
$d\co\pi\to sl(n,\CC)$ satisfying 
\[d(\gamma_1\gamma_2)=d(\gamma_1)+\Ad\circ\rho(\gamma_1)(d(\gamma_2))\quad,\ \forall\ \gamma_1,\ \gamma_2\in\pi\,.\]

It was observed by Andr\'e Weil \cite{Weil} that there is a natural inclusion of the Zariski tangent space 
$T_\rho^{Zar}(R_n(\pi))\hookrightarrow Z^1(\pi,sl(n,\CC)_\rho)$. Informally speaking, given a smooth curve $\rho_\epsilon$ of representations through $\rho_0=\rho$ one gets a $1$-cocycle $d\co\pi\to sl(n,\CC)$ by defining
\[ d(\gamma) := \left.\frac{d \, \rho_{\epsilon}(\gamma)}
{d\,\epsilon}\right|_{\epsilon=0} \rho(\gamma)^{-1},
\quad\forall\gamma\in\pi\,.\]

It is easy to see that the tangent space to the orbit by conjugation corresponds to the space of $1$-coboundaries $B^1(\pi,sl(n,\CC)_\rho)$. Here, $b\co\pi\to sl(n,\CC)$ is a coboundary if there exists 
$x\in sl(n,\CC)$ such that $b(\gamma)=\Ad\circ\rho(\gamma)(x)-x$. A detailed account can be found in \cite{L-M}.

Let $\dim_\rho R_n(\pi)$ be the local dimension of $R_n(\pi)$ at $\rho$ (i.e.\ the maximal dimension of the irreducible components of $R_n(\pi)$ containing $\rho$ \cite[Ch.~II]{Sh}). In the sequel we will use the following lemmas from \cite{H-P}:
\begin{lemma}\label{smoothness}
Let $\rho\in R_n(\pi)$ be given. If 
$\dim_\rho R_n(\pi)=\dim Z^1(\pi,sl(n,\CC)_\rho)$ then $\rho$ is a smooth point of the representation variety $R_n(\pi)$ and $\rho$ is contained in a unique component of $R_n(\pi)$ of dimension 
$\dim Z^1(\pi,sl(n,\CC)_\rho)$.
\end{lemma}

\begin{lemma}\label{leminj}
Let $A$ be a $\pi$-module and let $M$ be any $CW$-complex with $\pi_1(M)\cong\pi$. Then there are natural morphisms $H_i(M,A)\to H_i(\pi,A)$ which are isomorphisms for $i=0,1$ and a surjection for $i=2$. In cohomology there are natural morphisms $H^i(\pi,A)\to H^i(M,A)$ which are isomorphisms for $i=0,1$ and an injection for $i=2$.
\end{lemma}
\begin{remark}
Let $A$ be a $\pi$-module and $X$ a knot complement in $S^3$. The homomorphisms $H^*(\pi,A)\to H^*(X,A)$ and $H^*(\pi_1(\partial X),A)\to H^*(\partial X,A)$ are isomorphisms. This is a consequence of the asphericity of $X$ and $\partial X$. Moreover, the knot complement $X$ has the homotopy type of a $2$-dimensional CW-complex which implies that $H^k(\pi,A)=0$ for $k\geq3$.
\end{remark}

\section{Deforming metabelian representations}
\label{deform}

The aim of the following sections is to prove that, when $\alpha$ is a root of the Alexander polynomial of order $2$ then certain metabelian representations are smooth points of the representation variety.

In order to construct deformations of $\tilde{\rho}$ we use the classical approach, i.e. we first solve the corresponding formal problem and apply then a deep theorem of Artin \cite{A}. The formal deformations of a representation $\rho\co\pi\to SL(3,\CC)$ are in general determined by an infinite sequence of obstructions (see \cite{BA} and \cite{Gol}).

Given a cocycle in $Z^1(\pi, sl(n,\CC)_\rho)$ the first obstruction to integration is the cup product with itself. In general when the $k$-th obstruction vanishes, the obstruction of order $k+1$ is defined, it lives in $H^2(\pi, sl(n,\CC)_\rho)$.

Let $\rho\co\pi\to SL(n,\CC)$ be a representation. A formal deformation of $\rho$ is a homomorphism 
$\rho_\infty\co\pi\to SL(n,\CC[[t]])$

\[\rho_\infty(\gamma)=\exp\left(\displaystyle\sum_{i=1}^{\infty}t^iu_i(\gamma)\right)\rho(\gamma)\]
where $u_i\co\pi\to sl(n,\CC)$ are elements of $C^1(\pi,sl(n,\CC)_\rho)$ such that $ev_0\circ\rho_\infty=\rho$. 
Here $ev_0\co SL(n,\CC[[t]])\to SL(n,\CC)$ is the evaluation homomorphism at $t=0$ and $\CC[[t]]$ denotes the ring of formal power series. We will say that $\rho_\infty$ is a formal deformation up to  order $k$ of $\rho$ if $\rho_\infty$ is a homomorphism modulo $t^{k+1}$.

An easy calculation gives that $\rho_\infty$ is a homomorphism up to first order if and only if $u_1\in Z^1(\pi,sl(n,\CC)_\rho)$ is a cocycle. We call a cocycle $u_1\in Z^1(\pi,sl(n,\CC)_\rho)$ \emph{integrable} if there is a formal deformation of $\rho$ with leading term $u_1$.

\begin{lemma}
Let $u_1,\ldots,u_k\in C^1(\pi,sl(n,\CC)_\rho)$ such that
\[\rho_k(\gamma)=\exp\left(\displaystyle\sum_{i=1}^k t^iu_i(\gamma)\right)\rho(\gamma)\]
is a homomorphism into $SL(n,\CC[[t]])$ modulo $t^{k+1}$. Then there exists an obstruction class $\zeta_{k+1}:=\zeta_{k+1}^{(u_1,\ldots,u_k)}\in H^2(\pi, sl(n,\CC)_\rho)$ with the following properties:
\begin{description}
\item [(i)] There is a cochain $u_{k+1}\co\pi\to sl(n,\CC)_\rho$ such that
\[\rho_{k+1}(\gamma)=\exp\left(\displaystyle\sum_{i=1}^{k+1}t^iu_i(\gamma)\right)\rho(\gamma)\]
is a homomorphism modulo $t^{k+2}$ if and only if $\zeta_{k+1}=0$.
\item [(ii)] The obstruction $\zeta_{k+1}$ is natural, i.e. if $f\co\Gamma\to\pi$ is a homomorphism then $f^*\rho_k:=\rho_k\circ f$ is also a homomorphism modulo $t^{k+1}$ and $f^*(\zeta_{k+1}^{(u_1,\ldots,u_k)})=\zeta_{k+1}^{(f^*u_1,\ldots,f^*u_k)}$.
\end{description}

\end{lemma}

\begin{proof}
The proof is completely analogous to the proof of 
Proposition~3.1 in \cite{H-P-S}. We replace $SL(2,\CC)$ (resp. $sl(2,\CC)$) by 
$SL(n,\CC)$ (resp. $sl(n,\CC)$). 
\end{proof}

Let $i\co\partial X\to X$ be the inclusion. For the convenience of the reader, we state the following result which is implicitly contained in \cite{H-P}:

\begin{thm}\label{mainthm}
Let $\rho\in R_n(\pi)$ be a representation such that 
$H^0( X,sl(n,\CC)_\rho)=0$ i.e.\ the centralizer 
$Z(\rho)\subset SL(n,\CC)$ is finite.

If 
$\dim H^0(\partial X,sl(n,\CC)_\rho)=\dim H^2(X,sl(n,\CC)_\rho)=n-1$
and if $\rho\circ i_\#$ is a smooth point of $R_n(\pi_1(\partial X))$, then $\rho$  is a smooth point of the representation variety $R_n(\pi)$; it is contained in a unique irreducible component of dimension 
$n^2+n-2=(n+2)(n-1)$.
\end{thm}

\begin{proof}
Recall that the Zariski tangent space of $R_n(\pi)$ at $\rho$ is contained in $Z^1(\pi, sl(n,\CC)_\rho)$ \cite{Weil}. To prove the smoothness, we show that all cocycles in $Z^1(\pi, sl(n,\CC)_\rho)$ are integrable. Therefore, we prove that all obstructions vanish, by using the fact that the obstructions vanish on the boundary.

We consider the exact sequence in cohomology for the pair $(X,\partial X)$: 
\begin{multline*}
0\to H^0(\partial X,sl(n,\CC)_\rho)\to
H^1(X,\partial X,sl(n,\CC)_\rho)\to \\
 H^1(X,sl(n,\CC)_\rho)
\to H^1(\partial X,sl(n,\CC)_\rho)\to H^2(X,\partial X,sl(n,\CC)_\rho)\to\\
H^2(X,sl(n,\CC)_\rho)\stackrel{i_1^*}{\to}H^2(\partial X,sl(n,\CC)_\rho)\to 0\,.\end{multline*}
Poincar\'e duality implies that $\dim H^2(\partial X,sl(n,\CC)_\rho)=n-1$ and the Poincar\'e--Lefschetz duality  gives
\[ H^3(X,\partial X,sl(n,\CC)_\rho)\simeq H^0(X,sl(n,\CC)_\rho)^*=0\,.\]

Since $i_1^*$ is surjective and $\dim H^2(X,sl(n,\CC)_\rho)=\dim H^2(\partial X,sl(n,\CC)_\rho)$ we get $H^2(X,sl(n,\CC)_\rho)\simeq H^2(\partial X,sl(n,\CC)_\rho)$. From Lemma~\ref{leminj} we deduce that
\[i^*\co H^2(\pi,sl(n,\CC)_\rho)\to H^2(\pi_1(\partial X),sl(n,\CC)_\rho)\]
is an isomorphism.

We will now prove that every element of $Z^1(\pi, sl(n,\CC)_\rho)$ is integrable. Let $u_1,\ldots,u_k\co\pi\to sl(n,\CC)$ be given such that \[\rho_k(\gamma)=\exp\left(\displaystyle\sum_{i=1}^kt^iu_i(\gamma)\right)\rho(\gamma)\]
 is a homomorphism modulo $t^{k+1}$. Then the restriction $\rho_k\circ i_{\#}\co\pi_1(\partial X)\to SL(n,\CC[[t]])$ is also a formal deformation of order $k$.
Since $\rho\circ i_{\#}$ is a smooth point of the representation variety $R_n(\ZZ\oplus\ZZ)$, the formal implicit function theorem gives that $\rho_k\circ i_{\#}$ extends to a formal deformation of order $k+1$ (see \cite[Lemma~3.7]{H-P-S}). Therefore, we have that
\[0=\zeta_{k+1}^{(i^*u_1,\ldots,i^*u_k)}=i^*\zeta_{k+1}^{(u_1,\ldots,u_k)}\]
Now, $i^*$ is injective and the obstruction vanishes.

Hence all cocycles in $Z^1(\pi, sl(n,\CC)_\rho)$ are integrable. By applying Artin's theorem \cite{A} we obtain from a formal deformation of $\rho$ a convergent deformation (see \cite[Lemma~3.3]{H-P-S} or \cite[§~4.2]{BA}).

Thus $\rho$ is a smooth point of the representation variety $R_n(\pi)$. The Euler characteristic $\chi(X)$ vanishes. Hence, $\dim H^1(\pi,sl(n,\CC)_\rho)=\dim H^2(\pi,sl(n,\CC)_\rho)=n-1$. Since $\dim B^1(\pi,sl(n,\CC)_\rho)=n^2-1$, then $\dim Z^1(\pi,sl(n,\CC)_\rho)=n^2+n-2$.
\end{proof}

Following \cite[§~3.5]{Stein} we call an  $A\in SL(n,\CC)$ a
\emph{regular element} if the dimension of the centralizer  $Z(A)$ of $A$ in $SL(n,\CC)$ is $n-1$. Moreover note that $A$ is regular iff $Z(A)$ is abelian.
The regular elements form an open dense set in $SL(n,\CC)$
(see \cite[§~3.5]{Stein} for details).

\begin{lemma}\label{lembound}
Let $\rho\in R_n(\ZZ\oplus\ZZ)$ be a representation and let 
$\mu\in\ZZ\oplus\ZZ$ be simple i.e. there exists 
$\lambda\in\ZZ\oplus\ZZ$ such that $(\mu,\lambda)$ is a basis.

If  $\rho(\mu)\in SL(n,\CC)$ is a regular element then $\rho$ is a smooth point of 
$R_n(\ZZ\oplus\ZZ)$. 
It belongs to a $(n+2)(n-1)$-dimensional component of 
$R_n(\ZZ\oplus\ZZ)$.
\end{lemma}
\begin{proof}
We have that 
\[H^0(\ZZ\oplus\ZZ,sl(n,\CC)_\rho) = sl(n,\CC)^{\ZZ\oplus\ZZ}
= sl(n,\CC)^{\rho(\mu)}\cap sl(n,\CC)^{\rho(\lambda)}\,.\]
The regularity of $\rho(\mu)$ implies that 
$\dim sl(n,\CC)^{\rho(\mu)} = n-1$ and $\dim \CC[\rho(\mu)] =n$
where $\CC[\rho(\mu)]\subset M(n,\CC)$ denotes the algebra generated by $\rho(\mu)$ (see \cite[§~3.5]{Stein}). On the other hand we have
\[  \CC[\rho(\mu)]\cap sl(n,\CC) \subset sl(n,\CC)^{\rho(\mu)} \]
and the equality of the dimensions gives
$\CC[\rho(\mu)]\cap sl(n,\CC) = sl(n,\CC)^{\rho(\mu)}$.
Therefore we have for each $A\in Z(\rho(\mu))$ that
$ sl(n,\CC)^{\rho(\mu)}\subset sl(n,\CC)^A$ and 
\[ H^0(\ZZ\oplus\ZZ,sl(n,\CC)_\rho) = sl(n,\CC)^{\ZZ\oplus\ZZ}
= sl(n,\CC)^{\rho(\mu)}\cap sl(n,\CC)^{\rho(\lambda)} = sl(n,\CC)^{\rho(\mu)}\]
follows. This gives $\dim H^0(\ZZ\oplus\ZZ,sl(n,\CC)_\rho)=n-1$.
Now, Poincaré duality implies that 
$\dim H^2(\ZZ\oplus\ZZ,sl(n,\CC)_\rho)=n-1$ and since the Euler characteristic $\chi(\partial X)$ vanishes we obtain 
$\dim H^1(\ZZ\oplus\ZZ,sl(n,\CC)_\rho)=2(n-1)$. Thus 
$\dim B^1(\ZZ\oplus\ZZ,sl(n,\CC)_\rho)=n^2-n$ and 
$\dim Z^1(\ZZ\oplus\ZZ,sl(n,\CC)_\rho)=n^2+n-2$.

So to prove that $\rho$ is a smooth point of $R_n(\ZZ\oplus\ZZ)$, we will verify that $\rho$ is contained in a $(n^2+n-2)$-dimensional component of $R(\ZZ\oplus\ZZ)$.
Since the set of regular elements of $SL(n,\CC)$ is  open and since the dimension of the centralizer of a regular element is by definition $n-1$ it follows easily that $\rho$ is contained in an $(n^2-1) + (n-1)$ dimensional component.
\end{proof}

From now on, $sl(3,\CC)$ is considered as a $\pi$-module via the action of $\Ad\circ\tilde\rho$. Let's recall that the image of a meridian is given by \[\tilde\rho(\mu)=\alpha^{-1/3}\begin{pmatrix}
\alpha&0&0\\
0&1&1\\
0&0&1
\end{pmatrix}\,.\] 
which is a regular element of $SL(3,\CC)$.

We have:
\[H^0(\partial X, sl(3,\CC)) = sl(3,\CC)^{\tilde\rho(\mu)}\]
is two dimensional.

The next proposition will be proved in Section~\ref{cohometab}.
\begin{prop}\label{dimsl}
Let $K\subset S^3$ be a knot and suppose that the
$(t-\alpha)$-torsion of the Alexander module of $K$ is  of the form 
$\tau_\alpha=\CC[t,t^{-1}]\big/(t-\alpha)^2$.

Then we have:
\begin{enumerate}
\item $H^0(\pi,sl(3,\CC))=0\,.$
\item $\dim H^1(\pi,sl(3,\CC))=\dim H^2(\pi,sl(3,\CC))=2\,.$
\end{enumerate}
\end{prop}

\begin{proof}[First part of the proof of Theorem~\ref{mainresult}]
Recall that $\dim H^0(\partial X,sl(3,\CC))=2$. Now, we apply Theorem \ref{mainthm}, using Lemma \ref{lembound} and the fact that $\dim H^2(\pi,sl(3,\CC))=2$ (Proposition~\ref{dimsl}), to prove that $\tilde\rho$ is a smooth point of $R_3(\pi)$. 
By Proposition~\ref{dimsl} we have  
$\dim H^0(\pi,sl(3,\CC))=0$ and $\dim H^1(\pi,sl(3,\CC))=2$ which implies that $\dim Z^1(\pi,sl(3,\CC))=10$. Hence the representation 
$\tilde\rho$ is contained in a $10$-dimensional component $R_{\tilde\rho}$ of the representation variety (Lemma \ref{smoothness}). In Theorem~\ref{thm:stable} and Corollary~\ref{cor:nonmetab}, we will see that the component $R_{\tilde\rho}$ contains irreducible non metabelian representations.
\end{proof}

\section{Cohomology of metabelian representations}
\label{cohometab}

Throughout this section we will suppose that the
$(t-\alpha)$-torsion of the Alexander module of $K$ is  of the form 
$\tau_\alpha=\CC[t,t^{-1}]\big/(t-\alpha)^2$.
This implies that $\dim H^1(\pi,\CC_\alpha)=1$.

Recall that the Lie algebra $sl(3,\CC)$ is a $\pi$-module via the action of the adjoint representation $\Ad\circ\tilde\rho=\Ad\circ\rho_0$, where
\[\rho_0(\gamma)=\begin{pmatrix}
\alpha^{|\gamma|}&z(\gamma)&g(\gamma)\\
0&1&h(\gamma)\\
0&0&1
\end{pmatrix}\,.\]

The goal in this section is to prove Proposition~\ref{dimsl}.
The calculation of the dimensions of $H^*(\pi,sl(3,\CC))$
uses several  long exact sequences in cohomology associated to
the $\pi$-module $sl(3,\CC)$.

\subsection{The setup}
Denote by $(E_{ij})_{1\leq i,j\leq 3}$ the canonical basis of $M(3,\CC)$ and let $D_1=E_{22}-2E_{11}+E_{33}$ and 
$D_2=E_{11}-2E_{33}+E_{22}$. Then 
$(D_1, D_2, E_{ij}\mid 1\leq i\neq j\leq 3)$ is a basis of $sl(3,\CC)$.

It is not hard to check that: 
\[ \begin{array}{lllll}
H^0(\pi,\CC)\simeq\CC &;& H^1(\pi,\CC)\simeq\CC &;& H^2(\pi,\CC)=0\cr
H^0(\pi,\CC_\alpha)=0 &;& H^1(\pi,\CC_\alpha)\simeq\CC &;& H^2(\pi,\CC_\alpha)\simeq\CC
\end{array}
\]
For more details see \cite{BA-thesis}.

We define $\CC_+(3):=\left\langle E_{12},E_{23},E_{13}\right\rangle\ \text{and}\ b_+:=\left\langle D_1,D_2,E_{12},E_{23},E_{13}\right\rangle$ the Borel subalgebra of upper triangular matrices. Then the action of $\Ad\circ\rho_0$ on $\CC_+(3)$ and $b_+$ is given by
\begin{align}\label{actb}
\notag
\gamma\cdot E_{12}&=\alpha^{|\gamma|}E_{12}-\alpha^{|\gamma|}h(\gamma)E_{13}\\ \notag
\gamma\cdot E_{13}&=\alpha^{|\gamma|}E_{13}\\
\gamma\cdot E_{23}&=E_{23}+z(\gamma)E_{13}\\ \notag
\gamma\cdot D_1&=D_1+3z(\gamma)E_{12}-\left(3z(\gamma)h(\gamma)-3g(\gamma)\right)E_{13}\\
\gamma\cdot D_2&=D_2-3h(\gamma)E_{23}-3g(\gamma)E_{13} \notag
\end{align}
Hence $\CC_+(3)$ and $b_+$ are $\pi$-submodules of $sl(3,\CC)$. Moreover $\left\langle E_{13}\right\rangle$ is a $\pi$-submodule of $\CC_+(3)$ where the action is given by the multiplication by $\alpha$.  We have the following short exact sequence
\begin{equation}\label{suitecplus-old}
 0\to\left\langle E_{13}\right\rangle\to\CC_+(3)\to\CC_+(3)/\left\langle E_{13}\right\rangle\to 0\,.
\end{equation}
Note that $\langle E_{13}\rangle\simeq\CC_\alpha$ and that $\CC_+(3)/\left\langle E_{13}\right\rangle\simeq\CC\oplus\CC_\alpha$. The first isomorphism is induced by the  projection $p_{13}\co\langle E_{13}\rangle\to\CC_\alpha$  and the second isomorphism is induced by the projection 
\[ pr_1\co\CC_+(3)\to\CC\oplus\CC_\alpha,\quad
pr_1(M)=(p_{23}(M),p_{12}(M))\,.\]
Here $p_{ij}\co M(3,\CC)\to \CC$ denotes the projection onto the 
$(i,j)$-coordinates. Hence (\ref{suitecplus-old}) gives the short exact sequence
\begin{equation}\label{suitecplus}
 0\to \CC_\alpha \stackrel{i_1}{\longrightarrow}\CC_+(3)\stackrel{pr_1}{\longrightarrow}
 \CC\oplus\CC_\alpha\to 0
\end{equation}
where $i_1\co \CC_\alpha \to \CC_+(3)$ is given by $i_1(c)=c E_{13}$.

On the other hand, $\CC_+(3)$ is a $\pi$-submodule of $b_+$. Let us denote $D_+:=b_+/\CC_+(3)$, so we have the short exact sequence
 \begin{equation}\label{suiteb-old}
0\to\CC_+(3)\stackrel{i_2}{\longrightarrow} b_+\to D_+\to 0\,.
\end{equation}
The projection of an element of $b_+$ onto its coordinates on $D_1$ and $D_2$ induces by (\ref{actb}) an isomorphism 
$pr_2\co D_+\to\CC\oplus\CC$. Hence (\ref{suiteb-old}) gives the short exact sequence
\begin{equation}\label{suiteb}
0\to\CC_+(3)\stackrel{i_2}{\longrightarrow} b_+ \stackrel{pr_2}{\longrightarrow}\CC\oplus\CC \to 0\,.
\end{equation}

We define $\CC_-(3)$ as $sl(3,\CC)/b_+$.  Then we have a short exact sequence of $\pi$-modules
\begin{equation}\label{seqsl}
0\to b_+\to sl(3,\CC)\to\CC_-(3)\to 0\,.
\end{equation}
The action of $\Ad\circ\rho_0$ on the lower triangular matrices in 
$sl(3,\CC)$ is given by:
\begin{align}\label{actinf}
\notag
\gamma\cdot E_{21}&=\alpha^{-|\gamma|}E_{21}\pmod{b_+}\\
\gamma\cdot E_{31}&=\alpha^{-|\gamma|}h(\gamma)E_{21}+
\alpha^{-|\gamma|}E_{31}-\alpha^{-|\gamma|}z(\gamma)E_{32}\pmod{b_+}\\  \notag
\gamma\cdot E_{32}&=E_{32}\pmod{b_+}\,.
\end{align}
Let $\overline{E}_{ij}=E_{ij}+b_+$, for $1\leq j<i\leq 3$.
Equation~(\ref{actinf}) gives that 
$\langle \overline{E}_{21},\overline{E}_{32}\rangle$  is a 
$\pi$-submodule of $\CC_-(3)$ and that 
\[
\langle \overline{E}_{21},\overline{E}_{32}\rangle\simeq\CC_{\alpha^{-1}}\oplus\CC\,.\] Moreover, the quotient 
$\CC_-(3)/\langle \overline{E}_{21},\overline{E}_{32}\rangle$ is isomorphic to $\CC_{\alpha^{-1}}$. This isomorphism is simply induced by the projection $p_{31}$.
Hence we obtain a short exact sequence
\begin{equation}\label{suitecmoins}
0\to\CC_{\alpha^{-1}}\oplus\CC\to \CC_-(3)\to\CC_{\alpha^{-1}}\to 0\,.\end{equation}

\subsection{The computations}
\begin{lemma}\label{dimcplus}
Same assumptions as in Proposition~\ref{dimsl}.

For the cohomology groups $H^k(\pi,\CC_+(3))$ the following holds:
\[H^k(\pi,\CC_+(3))=0,\text{ if } k\neq 1,2,\]
\[(pr_1)_*\co H^1(\pi,\CC_+(3))\stackrel{\simeq}{\longrightarrow} H^1(\pi,\CC)\oplus H^1(\pi,\CC_\alpha)\]
is an isomorphism and there is a short exact sequence
\[0\to H^2(\pi,\CC_\alpha)\stackrel{(i_1)_*}{\longrightarrow} H^2(\pi,\CC_+(3))\stackrel{(pr_1)_*}{\longrightarrow} H^2(\pi,\CC_\alpha)\to 0\,.\]
In particular, $\dim H^1(\pi,\CC_+(3))=\dim H^2(\pi,\CC_+(3))=2$.
\end{lemma}

\begin{proof}
The long exact cohomology  sequence associated to (\ref{suitecplus}) gives:
\begin{displaymath}\begin{split}
0\rightarrow H^0(\pi,\CC_+(3))\rightarrow H^0(\pi,\CC)\stackrel{\delta^1}{\longrightarrow}H^1(\pi,\CC_\alpha)
\rightarrow H^1(\pi,\CC_+(3))\\
\rightarrow H^1(\pi,\CC)\oplus H^1(\pi,\CC_\alpha)\stackrel{\delta^2}{\rightarrow}H^2(\pi,\CC_\alpha)
\rightarrow H^2(\pi,\CC_+(3))\rightarrow  H^2(\pi,\CC_\alpha)\rightarrow 0\,.
\end{split}\end{displaymath}
In order to calculate $\delta^1\co H^0(\pi,\CC)\to H^1(\pi,\CC_\alpha)$
let $\delta$ denote the coboundary operator of $C^*(\pi,\CC_+(3))$.

If $c\in Z^0(\pi,\CC)=H^0(\pi,\CC)$ then 
$cE_{23}\in (pr_1)_*^{-1}(c,0) \subset C^0(\pi,\CC_+(3))$
and by (\ref{actb}) we obtain:
 \[ \delta(c E_{23})(\gamma)=c (\gamma-1)\cdot E_{23}= c z(\gamma)E_{13}\,.\]
Therefore $\delta^1(c)= c \{z\}\in H^1(\pi,\CC_\alpha)$.
Since $\{z\}\neq 0$ in $H^1(\pi,\CC_\alpha)$, $\delta^1$ is injective and hence an isomorphism (recall that $\dim H^0(\pi,\CC)=\dim H^1(\pi,\CC_\alpha)=1$). This implies that $H^0(\pi,\CC_+(3))=0$ and
the long exact sequence in cohomology becomes
\begin{equation}\label{suitec}
\begin{split}
0\rightarrow H^1(\pi,\CC_+(3))\rightarrow H^1(\pi,\CC)\oplus H^1(\pi,\CC_\alpha)\stackrel{\delta^2}{\longrightarrow}H^2(\pi,\CC_\alpha)\rightarrow \\H^2(\pi,\CC_+(3))\rightarrow H^2(\pi,\CC_\alpha)\rightarrow 0\,.
\end{split}
\end{equation}

Next we consider 
$\delta^2\co H^1(\pi,\CC)\oplus H^1(\pi,\CC_\alpha)\to H^2(\pi,\CC_\alpha)$. For $h'\in Z^1(\pi,\CC)=H^1(\pi,\CC)=\hom(\pi,\CC)$, we have that $h' E_{23}\in (pr_1)_*^{-1}(h',0)\subset C^1(\pi,\CC_+(3))$ 
and (\ref{actb}) gives: 
\begin{align*} 
\delta(h' E_{23})(\gamma_1,\gamma_2)=&
\gamma_1\cdot(h'(\gamma_2)E_{23})-h'(\gamma_1\gamma_2)E_{23}+h'(\gamma_1)E_{23})\\=& z\cup h'(\gamma_1,\gamma_2) E_{13}\,.
\end{align*}
Hence $\delta^2(h')=\{z\cup h'\}$ and a similar computation for  $z'\in Z^1(\pi,\CC_\alpha)$ gives
$\delta^2(\{z'\})=\{-h\cup z'\}$.

So $\delta^2(\{h'\}+\{z'\})=\{z\cup h'\}-\{h\cup z'\}$.
Since $\alpha$ is not a simple root of the Alexander polynomial it follows that
$\delta^2\equiv 0$ (see \cite[Theorem~3.2]{BA-L}). We obtain, from (\ref{suitec}) the following  exact sequences:
\[0\rightarrow H^1(\pi,\CC_+(3))\stackrel{(pr_1)_*}{\longrightarrow}H^1(\pi,\CC)\oplus H^1(\pi,\CC_\alpha)\rightarrow 0\]
\[0\rightarrow H^2(\pi,\CC_\alpha)\stackrel{(i_1)_*}{\longrightarrow} H^2(\pi,\CC_+(3))\stackrel{(pr_1)_*}{\longrightarrow} H^2(\pi,\CC_\alpha)\rightarrow 0\]
from which we deduce that $\dim H^1(\pi,\CC_+(3))=\dim H^2(\pi,\CC_+(3))=2$.
\end{proof}

In the following lemma, we will compute the dimensions of 
$H^*(\pi,b_+)$.

\begin{lemma}\label{dimbplus} Same assumptions as in Proposition~\ref{dimsl}.

We have $\dim H^0(\pi,b_+)=0$ and 
$\dim H^1(\pi,b_+)=\dim H^2(\pi,b_+)=1$. Moreover, we have
$\Ker (i_2)_* = \Ker(pr_1)_*$ where
\[ (i_2)_*\co H^2(\pi,\CC_+(3))\to H^2(\pi,b_+) \text{ and } 
(pr_1)_*\co H^2(\pi,\CC_+(3))\to H^2(\pi,\CC_\alpha) \,. \]

\end{lemma}

\begin{proof}
The short exact sequence (\ref{suiteb})
 gives the following long exact cohomology sequence
\[\begin{split}
0\to H^0(\pi,b_+)\to H^0(\pi,\CC)\oplus H^0(\pi,\CC)\stackrel{\delta^1}{\to}H^1(\pi,\CC_+(3))\to H^1(\pi,b_+)\\
\to H^1(\pi,\CC)\oplus H^1(\pi,\CC)\stackrel{\delta^2}{\to}H^2(\pi,\CC_+(3))\to H^2(\pi,b_+)\to 0\,.
\end{split}\]

A calculation similar to the one in the last proof gives that $\delta^1$ is injective. Thus $H^0(\pi,b_+)=0$ and 
\[0\to H^1(\pi,b_+)\to H^1(\pi,\CC)\oplus H^1(\pi,\CC)\stackrel{\delta^2}{\to}H^2(\pi,\CC_+(3))\to H^2(\pi,b_+)\to 0\]
is exact.

Now we are interested in $\delta^2\co H^1(\pi,\CC)\oplus H^1(\pi,\CC)\to H^2(\pi,\CC_+(3))$. The element $h' D_1\in C^1(\pi,b_+)$ projects via $(pr_2)_*$ onto $(h',0)\in  Z^1(\pi,\CC)\oplus Z^1(\pi,\CC)$. Moreover:
\begin{align*}
\delta(h' D_1)(\gamma_1,\gamma_2)=& 
\gamma_1\cdot(h'(\gamma_2)D_1)-h'(\gamma_1\gamma_2)D_1+h'(\gamma_1)D_1)\\
=&3z\cup h'(\gamma_1,\gamma_2)E_{12}-
3((zh-g)\cup h')(\gamma_1,\gamma_2)E_{13}\,,
\end{align*}
where $\delta$ denotes the coboundary operator of $C^*(\pi,b_+)$.
Here  $zh\in C^1(\pi,\CC_\alpha)$ is simply defined by $
zh (\gamma):= z(\gamma)h(\gamma)$, for $\gamma\in\pi$.
Similarly, $h'' D_2\in C^1(\pi,b_+)$ projects onto 
$(0,h'')\in  Z^1(\pi,\CC)\oplus Z^1(\pi,\CC)$ and  
\begin{align*}
\delta(h'' D_2)(\gamma_1,\gamma_2)=&
\gamma_1\cdot(h''(\gamma_2)D_2)-h''(\gamma_1\gamma_2)D_2+
h''(\gamma_1)D_2)\\
=&-3h\cup h''(\gamma_1,\gamma_2)E_{23}-
3g\cup h''(\gamma_1,\gamma_2)E_{13}\,.
\end{align*}

We know that $\{h\cup h''\}=0$ since $H^2(\pi,\CC)=0$. 
So let $h_2\co\pi\to\CC$ be a $1$-cochain such that 
$\delta h_2+h\cup h''=0$. 
Then $h'' D_2 + 3h_2 E_{23}\in C^1(\pi,b_+)$ projects also via 
$(pr_2)_*$ onto $(0,h'')\in  Z^1(\pi,\CC)\oplus Z^1(\pi,\CC)$ and
\[\delta(h'' D_2 + 3h_2 E_{23})(\gamma_1,\gamma_2)=
3(z\cup h_2(\gamma_1,\gamma_2)+g\cup h''(\gamma_1,\gamma_2))E_{13}\,.\]
Hence 
\[\delta^2(0,h'')=3 \{z\cup h_2+g\cup h''\}E_{13} \in 
(i_1)_*\big(H^2(\pi,\CC_\alpha)\big)\subset H^2(\pi,\CC_+(3))\,.\]
Moreover, we know that $\{z\cup h_2+g\cup h''\}\neq 0$ in 
$H^2(\pi,\CC_\alpha)$ (see \cite[Theorem~1]{Hajer} which implies $\rk \delta^2\geq 1$.

Similarly, there exists $g'\co\pi\to\CC_\alpha$ a $1$-cochain satisfying $\delta g'+z\cup h'=0$ and $h' D_1+ 3g' E_{12}\in C^1(\pi,b_+)$ projects also onto $(h',0)\in  Z^1(\pi,\CC)\oplus Z^1(\pi,\CC)$.
We obtain:
\[\delta(h' D_1+ 3g' E_{12})(\gamma_1,\gamma_2)=
-3(h\cup g'(\gamma_1,\gamma_2)+
(zh-g)\cup h'(\gamma_1,\gamma_2))E_{13}\,.\]
Note that $\delta(zh-g) +h\cup z =0$.

This gives
$\image\delta^2\subset (i_1)_*(H^2(\pi,\CC_\alpha))\subset 
H^2(\pi,\CC_+(3))$. In particular, $\rk \delta^2\leq 1$ and hence 
$\rk\delta^2=1$. Moreover we have $\image \delta^2 = \image(i_1)_*$ and hence $\Ker (i_2)_* = \Ker (pr_1)_*$.

The long exact sequence in cohomology gives 
$\dim H^1(\pi,b_+)=\dim H^2(\pi,b_+)=1$.
\end{proof}

\begin{lemma}\label{dimcmoins}
The short exact sequence (\ref{suitecmoins}) implies that
$H^0(\pi,\CC_-(3))\simeq H^0(\pi,\CC)$ and gives the following exact sequences:
\[0\to H^1(\pi,\CC_{\alpha^{-1}})\oplus H^1(\pi,\CC)\to H^1(\pi,\CC_-(3))\to H^1(\pi,\CC_{\alpha^{-1}})\to 0\]
\[0\to H^2(\pi,\CC_{\alpha^{-1}})\to H^2(\pi,\CC_-(3))\to H^2(\pi,\CC_{\alpha^{-1}})\to 0\,.\]
In particular, we have $\dim H^0(\pi,\CC_-(3))=1$, 
$\dim H^1(\pi,\CC_-(3))=3$ and $\dim H^2(\pi,\CC_-(3))=2$.

\end{lemma}

\begin{proof}
The long exact cohomology sequence associated to the short exact sequence (\ref{suitecmoins}) and the fact that 
$H^0(\pi,\CC_{\alpha^\pm})=0$ gives the isomorphism
\[0\to H^0(\pi,\CC)\stackrel{\simeq}{\longrightarrow}H^0(\pi,\CC_-(3))\to 0\]
and the exact sequence
\begin{multline*}
0\stackrel{\delta^1}{\longrightarrow}H^1(\pi,\CC)\oplus H^1(\pi,\CC_{\alpha^{-1}})\to H^1(\pi,\CC_-(3))\to \\ H^1(\pi,\CC_{\alpha^{-1}})\stackrel{\delta^2}{\longrightarrow}H^2(\pi,\CC_{\alpha^{-1}})
\to H^2(\pi,\CC_-(3))\to H^2(\pi,\CC_{\alpha^{-1}})\to 0\,.\end{multline*}

Now, by similar computation as before and by using the fact that
$H^2(\pi,\CC)=0$ we obtain
\[\delta^2(\{z_-\})=\{ h\cup z_-\}\in H^2(\pi,\CC_{\alpha^{-1}})\]
Since $\alpha$ is a double root of the Alexander polynomial, $\{h\cup z_-\}=0$ (see \cite {H-P} or \cite{Hajer}), so $\delta^2\equiv 0$ and lemma follows.

\end{proof}

\begin{proof}[Proof of Proposition~\ref{dimsl}]
The short exact sequence~(\ref{seqsl}) of $\pi$-modules gives the following long exact cohomology sequence
\begin{multline*}
0\to H^0(\pi,sl(3,\CC))\to H^0(\pi,\CC_-(3))\stackrel{\delta^1}{\longrightarrow}H^1(\pi,b_+)\to H^1(\pi,sl(3,\CC))\\
\to H^1(\pi,\CC_-(3))\stackrel{\delta^2}{\longrightarrow}H^2(\pi,b_+)\to H^2(\pi,sl(3,\CC))\to H^2(\pi,\CC_-(3))\to 0\,.
\end{multline*}

An explicit calculation gives:
\[
H^0(\pi,sl(3,\CC))=\{A\in sl(3,\CC)\ | \ \gamma\cdot A=A,\ \forall\ \gamma\in\pi\}=\{0\}\]
which implies that $\delta^1$ is injective. Since $\dim H^0(\pi,\CC_-(3))=\dim H^1(\pi,b_+)$ (Lemmas~\ref{dimbplus} and~\ref{dimcmoins}), $\delta^1$ is an isomorphism. So, we obtain
\[0\to H^0(\pi,\CC_-(3))\stackrel{\simeq}{\longrightarrow}H^1(\pi,b_+)\to 0\] and
\begin{multline*}
0\to H^1(\pi,sl(3,\CC))\to H^1(\pi,\CC_-(3))\stackrel{\delta^2}{\to}H^2(\pi,b_+)\to\\ H^2(\pi,sl(3,\CC))\to H^2(\pi,\CC_-(3))\to 0\,.
\end{multline*}

Now, $H^1(\pi,\CC)\oplus H^1(\pi,\CC_{\alpha^{-1}})$ injects in $H^1(\pi,\CC_-(3))$ (Lemma \ref{dimcmoins}), so to understand the map $\delta^2$, we do the following calculations:
\[\gamma\cdot E_{32}=E_{32}+\frac{1}{3} h(\gamma)(2D_2+D_1)+g(\gamma)E_{12}-h^2(\gamma)E_{23}-g(\gamma)h(\gamma)E_{13}\,.\]
Hence for $h'\in Z^1(\pi,\CC)\simeq\hom(\pi,\CC)$, we have
$h'\overline{E}_{32} \in Z^1(\pi,\CC_-(3))$ and 
$h' E_{32}\in C^1(\pi, sl(3,\CC))$ projects onto $h'\overline{E}_{32}$.
Moreover, 
\[
\delta(h' E_{32})=
\frac{1}{3} h\cup h' (2D_2+D_1)+g\cup h' E_{12}-h^2\cup h' E_{23}-gh\cup h' E_{13}\,.\]
Here we let $\delta$ denote the coboundary operator of $C^*(\pi, sl(3,\CC))$.

Let $h_2\co\pi\to\CC$ be a $1$-cochain such that $\delta h_2+h\cup h'=0$, then 
\[ h' E_{32} +\frac 1 3 h_2 (2D_2+D_1)\in C^1(\pi, sl(3,\CC))\] 
projects also onto $h'\overline{E}_{32}$ and 
\begin{multline*}
\delta(h' E_{32} +\frac 1 3 h_2 (2D_2+D_1)) (\gamma_1,\gamma_2)=
(g\cup h'+z\cup h_2)(\gamma_1,\gamma_2)E_{12}\\-(h^2\cup h'+2h\cup h_2)(\gamma_1,\gamma_2)E_{23}
-(gh\cup h'+(zh+g)\cup h_2)(\gamma_1,\gamma_2)E_{13}\,.
\end{multline*}
This gives that 
\[\delta(h' E_{32} +\frac 1 3 h_2 (2D_2+D_1))\in
\image\big((i_2)_*\co H^2(\pi,\CC_+(3))\to H^2(\pi,b_+)\big)\,.\]
Moreover,
\[ (pr_1)_* \big( \delta(h' E_{32} +\frac 1 3 h_2 (2D_2+D_1)) \big)
= g\cup h'+z\cup h_2\,.\]
Since $\{g\cup h'+z\cup h_2\}\neq 0$ in $H^2(\pi,\CC_\alpha)$ 
(see \cite[Theorem~1]{Hajer}) we obtain from
$\Ker (pr_1)_*= \Ker (i_2)_*$ that
$\delta^2(h'\overline{E}_{32}) \neq 0$ and hence $\rk\delta^2\geq 1$. 
Moreover, $\dim H^2(\pi,b_+)=1$, so $\rk\delta^2=1$ and the long exact sequence enables us to conclude that
\[0\to H^1(\pi,sl(3,\CC))\to H^1(\pi,\CC_-(3)) \to H^2(\pi,b_+)\to 0\]
is exact and that
\[ H^2(\pi,sl(3,\CC))\stackrel{\simeq}{\longrightarrow} H^2(\pi,\CC_-(3))\]
is an isomorphism.
In particular, using Lemmas \ref{dimcplus}, \ref{dimbplus} and \ref{dimcmoins} we have 
$\dim H^1(\pi,sl(3,\CC))=\dim H^2(\pi,sl(3,\CC))=2$.
\end{proof}

\begin{remark}
If we consider the exact sequence in cohomology for the pair $(X,\partial X)$, we have:
\[H^1(X,\partial X,sl(3,\CC))\to H^1(X,sl(3,\CC))\stackrel{i_2^*}{\to}H^1(\partial X,sl(3,\CC))\,.\]
Applying the Poincar\'e duality, we obtain $\rk i_2^*=\frac{1}{2}\dim H^1(\partial X,sl(3,\CC))=2$. So $\dim H^1(X,sl(3,\CC))=\dim H^1(\pi,sl(3,\CC))\geq 2$.
\end{remark}

\section{The nature of the deformations}
\label{nature}

Throughout this section we will suppose that the
$(t-\alpha)$-torsion of the Alexander module of $K$ is  of the form 
$\tau_\alpha=\CC[t,t^{-1}]\big/(t-\alpha)^2$.

A representation $\rho\co\pi\to SL(n,\CC)$ is called \emph{reducible}
if there exists a proper subspace $V\subset \CC^n$ such that $\rho(\pi)$ preserves $V$. Otherwise $\rho$ is called \emph{irreducible}.
By Burnsides theorem, a representation $\rho$ is irreducible if and only if the image $\rho(\pi)$ generates the full matrix algebra $M(n,\CC)$. The \emph{orbit} of a representation $\rho$ is the  subset
$\mathcal{O} (\rho)=\{ \Ad_A\circ\rho\mid A\in SL(n,\CC)\}\subset R_n(\pi)$.

Note that the orbit of an irreducible representation is closed.
The orbit of the representation $\tilde \rho$ is not closed. This might be seen by looking at the one parameter subgroup $\lambda\co\CC^*\to SL(3,\CC)$ given by $\lambda (t)=\mathrm{diag}(t,1,t^{-1})$. It follows that
\[ \rho_\alpha(\gamma) := \lim_{t\to0} \lambda(t) \tilde\rho(\gamma)\lambda(t)^{-1}\]
is a diagonal representation $\rho_\alpha\co\pi\to SL(3,\CC)$  given by $\rho_\alpha(\mu) = \alpha^{-1/3} \mathrm{diag}(\alpha,1,1)$. 
Note that the orbit  $\mathcal O (\rho_\alpha)$ is closed and $4$-dimensional. It is contained in the closure 
$\overline{\mathcal O (\tilde \rho)}$ which is $8$-dimensional.

\begin{definition}
A representation $\rho\in R_n(\pi)$ is called \emph{stable} if its orbit 
$\mathcal{O} (\rho)$ is closed and if the isotropy group $Z(\rho)$ is finite. We denote by $S_n(\pi)\subset R_n(\pi)$ the set of stable representations.
\end{definition}

\begin{rem}
By a result of Newstead \cite[Proposition~3.8]{New78}, the set $S_n(\pi)$ is Zariski open in $R_n(\pi)$. However, $S_n (\pi)$ might be empty. 
\end{rem}

Next we will see that there are stable deformations of $\tilde\rho$. In order to proceed we will assume that there is a Wirtinger generator $S_1$  of $\pi$ such that $z(S_1)=0= g(S_1)$.
This can always be arranged by adding a coboundary to $z$ and $g$ i.e.\ by conjugating the representation $\tilde\rho$.

\begin{lemma} \label{lem:eigen}
Let $\rho_t\co\pi\to SL(3,\CC)$ be a curve in $R(\pi)$ with $\rho_0=\tilde{\rho}$. Then there exists a curve $C_t$ in $SL(3,\CC)$ such that $C_0 =I_3$ and 
\[ \Ad_{C_t}\circ \rho_t (S_1) = 
\begin{pmatrix}
a_{11}(t) & 0 & 0 \\
0 & a_{22}(t) & a_{23}(t)  \\
0 & a_{32}(t) & a_{33}(t)  
\end{pmatrix}\,.\]
for all sufficiently small $t$.
\end{lemma}
\begin{proof}
Let $A(t) := \rho_t (S_1) $.
Note that $\alpha^{2/3}$ is a simple root of the characteristic polynomial of $A(0)$. Hence there is a simple eigenvalue
$a_{11}(t)$ of $A(t)$ which depends analytically on $t$. Note that  the corresponding eigenvector $v_t$ can be chosen to depend also analytically on $t$ and such that $v_0$ is the first canonical basis vector $e_1$ of $\CC^3$. Hence $(v_t, e_2, e_3)$ is a basis for all sufficiently small $t$ and $A(t)$ takes the form
\[ A(t) = \begin{pmatrix}
a_{11}(t) & a_{12}(t) & a_{13}(t) \\
0 & a_{22}(t) & a_{23}(t)  \\
0 & a_{32}(t) & a_{33}(t)  
\end{pmatrix}\,.\]
Next observe that the matrix $(A_{11}(t) - a_{11}(t) I_2)$ is invertible for sufficiently small $t$. Here $A_{11}$ denotes the minor obtained from $A$ by eliminating the first row and the first column.
Hence the system
\[ (a_{12}(t),a_{13}(t))+(x(t),y(t))(A_{11}(t)-a_{11}(t)I_2)=0 \]
has a unique solution and for
\[ P(t) =
\begin{pmatrix}
1 &  x(t) & y(t) \\
0 & 1 & 0  \\
0 & 0 & 1
\end{pmatrix} \]
the matrix $P(t) A(t) P(t)^{-1}$ has the desired form.
\end{proof}

For the next step we choose a second Wirtinger generator $S_2$ of $\pi$ such that $z(S_2)=b\neq 0 =z(S_1)$. This is always possible since $z$ is not a coboundary. Hence
\[ \tilde\rho(S_1) = 
\begin{pmatrix}
\alpha^{2/3} & 0 & 0 \\
0 & \alpha^{-1/3} & \alpha^{-1/3}  \\
0 & 0 & \alpha^{-1/3}  
\end{pmatrix} \text{ and }
\tilde\rho(S_2) = 
\begin{pmatrix}
\alpha^{2/3} & b & c \\
0 & \alpha^{-1/3} & \alpha^{-1/3}  \\
0 & 0 & \alpha^{-1/3}  
\end{pmatrix}\]
where $b\neq 0$.

\begin{prop}\label{prop:irred}
Let  $A(t)$ and $B(t)=(b_{ij}(t))_{1\leq i,j\leq 3}$ be matrices depending analytically on $t$ such that 
\[A(t)=\begin{pmatrix}
a_{11}(t) & 0 & 0 \\
0 & a_{22}(t) & a_{23}(t)  \\
0 & a_{32}(t) & a_{33}(t)  
\end{pmatrix}, \quad
A(0) = \tilde\rho(S_1) \text{ and }
B(0) = \tilde\rho(S_2)\,.\]

If the first derivative $b'_{31}(0)\neq 0$ then for sufficiently small $t$, $t\neq 0$, the matrices $A(t)$ and $B(t)$ generate the full matrix algebra.

\end{prop}
\begin{proof} We denote by $\mathcal A_t\subset M(3,\CC)$ the algebra generated by $A(t)$ and $B(t)$.
Let $\chi_{A_{11}}(X)$ denote the characteristic polynomial of $A_{11}(t)$. It follows that $\chi_{A_{11}}(a_{11}(t))\neq 0$ for small $t$ and hence
\[ \frac{\chi_{A_{11}}(A(t))}{\chi_{A_{11}(t)}(a_{11}(t))}
=
\begin{pmatrix}
1 & 0 & 0 \\
0 & 0 & 0  \\
0 & 0 & 0  
\end{pmatrix}=
\begin{pmatrix}
1  \\
0   \\
0   
\end{pmatrix}\otimes
\begin{pmatrix}
1 & 0 & 0 \\
\end{pmatrix}
\in \CC[A(t)]\subset\mathcal A_t\,.\]

In the next step we will prove that
\[ \mathcal A_t \begin{pmatrix}
1  \\
0   \\
0   
\end{pmatrix} = \CC^3 \text{ and }
\begin{pmatrix}
1 & 0 & 0 \\
\end{pmatrix} \mathcal A_t =\CC^3\quad,\ \text{for small}\ t\in\CC^*\,.\]
It follows from this that $\mathcal A_t$ contains all rank one matrices 
since a rank one matrix can be written as $v\otimes w$ where $v$ is a column vector and $w$ is a row vector. Note also that
$ A (v \otimes w) = (A v) w$ and $ (v \otimes w) A = v\otimes (wA)$. Since each matrix is the sum of rank one matrices the proposition follows.

The vectors $\begin{pmatrix}
1 & 0 & 0 \\
\end{pmatrix} A(0)$, $\begin{pmatrix}
1 & 0 & 0 \\
\end{pmatrix} B(0)$ and $\begin{pmatrix}
1 & 0 & 0 \\
\end{pmatrix} B(0)^2$ form a basis of the space of row vectors. Hence
$\begin{pmatrix}
1 & 0 & 0 \\
\end{pmatrix} \mathcal A_t$ is the space of row vectors for sufficiently small $t$.

Consider the three column vectors
\[ a(t) := A(t) \begin{pmatrix}
1  \\
0   \\
0   
\end{pmatrix},\quad b(t):= B(t) \begin{pmatrix}
1  \\
0   \\
0   
\end{pmatrix} \text { and }
c(t) := A(t) b(t) \]
and define the function $f(t) := \det (a(t), b(t),c(t))$. It follows that
$f(t) = a_{11}(t) g(t) $ where $g(t)$ is given by
\[ g(t) = \begin{vmatrix}
b_{21}(t) & a_{22}(t) b_{21}(t) + a_{23}(t) b_{31}(t)\\
b_{31}(t) & a_{32}(t) b_{21}(t) + a_{33}(t) b_{31}(t)
\end{vmatrix}\,.\]
Now it is easy to see that $g(0)=g'(0)=0$ but $g''(0)= -\alpha^{-1/3} (b'_{31}(0))^2$. Hence $g(t)\neq0$ for sufficiently small $t$, $t\neq 0$.

\end{proof}

\begin{lemma}\label{lem:galois}
Let $z_{\pm}\in Z^1(\pi,\CC_{\alpha^{\pm 1}})$ be nontrivial cocycles such that $z_+(S_1)=z_-(S_1)=0$. If $z_+(S_2)\neq 0$ then $z_-(S_2)\neq 0$.
\end{lemma}
\begin{proof}
We define $a:= \alpha+\alpha^{-1}$. The 
number $a$ is defined over $\QQ$ since the Alexander polynomial is 
symmetric. 
Now we have an extension of degree two $\QQ(a) \subset \QQ(\alpha)$.
The defining equation is simply $x^{2} - a x + 1=0$ and we obtain
a Galois automorphism $\tau\co\QQ(\alpha)\to\QQ(\alpha)$ of order two with fixed field $\QQ(a)$ and $\tau(\alpha)=\alpha^{-1}$.

By fixing a Wirtinger presentation 
$\pi=\langle S_1,\ldots, S_n\mid R_1,\ldots,R_{n-1}\rangle$, each cocycle $z_\pm$ corresponds to a solution of a linear system 
$J(\alpha^{\pm 1}) \mathbf{z}=0$ where 
$J\in M_{n-1,n}(\ZZ[t,t^{-1}])$ is the Jacobi matrix of the presentation (see \cite[p.~976]{H-P}). More precisely, if $\mathbf{z} =(s_1,\ldots,s_n)$ is a solution of the system $J(\alpha^{\pm 1}) \mathbf{z}=0$ then the corresponding cocycle is given by
$z_\pm(S_i) = s_i$, for $1\leq i\leq n$.

If $\mathbf{z}_+$ is a solution of $J(\alpha) \mathbf{z}=0$ with 
$s_1=0$ and $s_2\neq 0$ then
$ \tau( \mathbf{z}_+)$ is a solution of  
$J(\alpha^{-1}) \mathbf{z}=0$ since $\tau$ is an automorphism it follows $\tau(s_1)=0$ and $\tau(s_2)\neq0$. Let
$\tilde z_-$ denote the cocycle given by $\tilde z_- (S_i) = \tau(s_i)$.
Note that $\tilde z_-$ is nontrivial since 
$\tilde z_-(S_1)\neq \tilde z_-(S_2)$.

It follows from Blanchfield duality that $\alpha^{-1}$ is a double root of the Alexander polynomial and that 
$\dim H^1(\pi,\CC_{\alpha^{-1}}) =1$ (see \cite[Proposition~4.7]{H-P}).
Hence if $ z_-$ is any nontrivial cocycle then there exists $t\in\CC^*$ and $b\in\CC$ such that 
$ z_-(S_i) = t \tilde z_-(S_i) + (\alpha^{-1} -1)b$. Now $z_-(S_1)=0$ implies that $b=0$ and hence $z_-(S_2) = t \tilde z_-(S_2)\neq 0$.
\end{proof}

Recall from the proof of Proposition~\ref{dimsl} that the projection 
$sl(3,\CC)\to sl(3,\CC)/b_+\cong\CC_-(3)$
induces an isomorphism 
\[ \Phi\co H^1(\pi,sl(3,\CC)) \stackrel{\cong}\to \Ker\big( H^1(\pi,\CC_-(3))\stackrel{\delta^2}{\longrightarrow}H^2(\pi,b_+)\big)\,.\]

Moreover, recall from Lemma~\ref{dimcmoins} that there is a short exact sequence
\[0\to H^1(\pi,\CC_{\alpha^{-1}})\oplus H^1(\pi,\CC)\to H^1(\pi,\CC_-(3))\to H^1(\pi,\CC_{\alpha^{-1}})\to 0\,.\]
In the sequel we will fix a non trivial cocycle 
$z_-\in Z^1(\pi,\CC_{\alpha^{-1}})$ such that $z_-(S_1)=0$. It follows from the preceding lemma that $z_-(S_2)\neq0$.
Moreover we have that the two cocycles 
$h\cup z_- \in Z^2(\pi,\CC_{\alpha^{-1}})$ and 
$z\cup z_- \in Z^2(\pi,\CC)$ are coboundaries. We will also fix cochains $g_-\co\pi\to\CC_{\alpha^{-1}}$ and $g_0\co\pi\to\CC$ such that 
\[ \delta g_- + h\cup z_- =0 \text{ and } \delta g_0 + z\cup z_- =0\,.\]
From Equation~(\ref{actinf}) and the above sequence we obtain that $H^1(\pi,\CC_-(3))$ is a three dimensional  vector space 
with basis 
\[ \bar z_1 = z_-\overline{E}_{21},\quad \bar z_2 = h \overline{E}_{32} \text{ and }
\bar z_3 = z_-\overline{E}_{31} - g_0 \overline{E}_{32} + 
g_- \overline{E}_{21}\,.\]
Hence every $z\in Z^1(\pi,sl(3,\CC))$ has the form
\[ z  =
\begin{pmatrix}
* & * & *\\
t_1  z_- + t_3 g_- + \delta b_1 & * & *\\
t_3  z_- + \delta b_3 & t_2  h - t_3 g_0  & *
\end{pmatrix}\]
where $t_i\in\CC$.

\begin{thm}\label{thm:stable}
There exist  deformations  
$\rho_t\co\pi\to SL(3,\CC)$ such that $\rho_0=\tilde{\rho}$, with the property  that $\rho_t$ is stable for all sufficiently small $t$, $t\neq 0$.
\end{thm}
\begin{proof}
Note that  $Z(\rho)$ is finite if and only if $H^0(\pi, sl(3,\CC)_\rho ) =0$. Moreover, the condition 
$H^0(\pi, sl(3,\CC)_\rho ) =0$ is an open condition on the representation variety. Hence all representation sufficiently close to $\tilde\rho$ have finite stabilizer.

Let  $z\in Z^1(\pi, sl(3,\CC))$ be a cocycle such that
$\Phi(z) = t_1 \bar z_1 +t_2 \bar z_2 + t_3 \bar z_3$ with $t_3\neq 0$. Such a cocycle exists always since 
$\bar z_2 \not\in \Ker \delta^2$ and $\dim \Ker \delta^2 = 2$.

Let $\rho_t$ be a deformation of $\tilde{\rho}$ with leading term 
$z$. We apply Lemma~\ref{lem:eigen} to this deformation for $A(t) = \rho_t(S_1)$ and $B(t)=\rho_t(S_2)$. 
Since $a_{31}(t)\equiv0$ it follows that  
\[ a'_{31}(0) = \alpha^{2/3} (t_3 z_-(S_1) + (\alpha^{-1}-1)b_3)=0 \] 
and hence $b_3 =0$. By Lemma~\ref{lem:galois} we obtain 
$b'_{31}(0) = \alpha^{2/3} t_3 z_-(S_2)\neq0$.

Hence we can apply Proposition~\ref{prop:irred} and obtain that
$\rho_t$ is irreducible for sufficiently small $t\neq0$.
\end{proof}

\begin{cor}\label{cor:nonmetab}
There exist irreducible, non metabelian deformations of $\tilde \rho$.
\end{cor}

\begin{proof}
Let $\rho_t$ be a deformation of $\tilde\rho$ such that $\rho_t$ is irreducible.
Then for sufficiently small $t$ we have that $\tr \rho_t (\mu)$ is close to $\tr \tilde\rho(\mu) = \alpha^{-1/3}(\alpha +2)$. Moreover we have 
$\tr \tilde\rho(\mu)\neq 0$ since $-2$ is not a root of the Alexander polynomial: $(x+2)\mid \Delta_K(x)$ implies $3\mid\Delta_K(1)=\pm1$ which is impossible.

By Theorem~1.2 of \cite{BF08}, we have for every irreducible metabelian representation $\rho\co\pi\to SL(3,\CC)$ that $\tr \rho(\mu)=0$. Hence $\rho_t$ is irreducible non metabelian for sufficiently small $t$.
\end{proof}

\begin{rem}
Let $\rho_\alpha\co\pi\to SL(3,\CC)$ be the diagonal representation given by $\rho_\alpha(\mu) = \alpha^{-1/3} \mathrm{diag}(\alpha,1,1)$. The orbit $\mathcal O (\rho_\alpha)$ is contained in the closure 
$\overline{\mathcal O (\tilde \rho)}$. Hence $\tilde\rho$ and $\rho_\alpha$ project to the same point $\chi_\alpha$ of the variety of characters
$X_3(\pi)= R_3(\pi)\sslash SL(3,\CC)$.

It would be natural to study the local picture of the variety of characters
$X_3(\pi)= R_3(\pi)\sslash SL(3,\CC)$ at $\chi_\alpha$ as done in \cite[§\ 8]{H-P}. Unfortunately, there are much more technical difficulties since in this case the quadratic cone $Q(\rho_\alpha)$ coincides with the Zariski tangent space $Z^1(\pi, sl(3,\CC)_{\rho_\alpha})$. Therefore the third obstruction has to be considered. 
\end{rem}

\end{document}